\documentclass[11pt]{amsart}
\usepackage{geometry,amssymb}                
\title[An inverse approach to hyperspheres of prescribed mean curvature]{An inverse approach to hyperspheres of prescribed mean curvature in Euclidean space}
\author{Paolo Caldiroli}
\address{Dipartimento di Matematica, Universit\`a di Torino, via Carlo Alberto, 10 -- 10123 Torino, Italy.}
\email{paolo.caldiroli@unito.it}
\thanks{The author is member of the Gruppo Nazionale per l'Analisi Matematica, la Probabilit\`a e le loro Applicazioni (GNAMPA) of the Istituto Nazionale di Alta Matematica (INdAM)}

\catcode`\@=11
\@addtoreset{equation}{section}

\catcode`\@=12

\newtheorem{Theorem}{Theorem}[section]
\newtheorem{Lemma}[Theorem]{Lemma}

\newtheorem{Remark}[Theorem]{Remark}
\newtheorem{Example}[Theorem]{Example}

\def\Proof{\noindent\textit{Proof. }}
\def\qed{$~\square$\goodbreak \medskip}

\begin{document}
\begin{abstract}
We construct families of smooth functions $H\colon\mathbb{R}^{n+1}\to\mathbb{R}$ such that the Euclidean $(n+1)$-space is completely filled by not necessarily round hyperspheres of mean curvature $H$ at every point.
\end{abstract}

\maketitle

\section*{Introduction}

This article deals with the problem of $n$-dimensional hypersurfaces of prescribed mean curvature in Euclidean space. More precisely, considering a function $H\colon\mathbb{R}^{n+1}\to\mathbb{R}$, one looks for compact, embedded hypersurfaces whose mean curvature at any point $p$ is given by $H(p)$. We focus our attention to hypersurfaces which are diffeomorphic to the sphere $\mathbb{S}^{n}$ and we call them \emph{$H$-bubbles}. 

When $H$ is a nonzero constant, $H$-bubbles are exactly round spheres of radius $|H|^{-1}$ centered at any point of the space. The case of nonconstant functions $H$ is much more delicate and has been intensively investigated in the last years. In fact the behaviour of $H$ crucially affects the existence  of $H$-bubbles (see, e.g., \cite{CaGu2007, 
CaMu2004, 
CaMu2011, KiLa2010, Mu2011, TrWe1983}). Moreover, concerning the location of $H$-bubbles, some constraints occur. For example, if $H$ depends just on one variable in a monotone way in some domain $D$ of the space, then no $H$-bubble can be found in $D$ (see \cite[Proposition 4.1]{CaMu2004}). In a perturbation setup, one sees that $H$-bubbles concentrate only at critical points of $H$ or of other functions related to $H$ (see, e.g., \cite[Proposition 4.3]{CaMu2004}, \cite[Theorem 2.1]{Ca2004}). 

In this paper, we are interested in a sort of inverse problem, more precisely, the existence of nonconstant curvature functions $H$ such that the whole Euclidean $(n+1)$-space can be filled by $H$-bubbles, in other words, for every point $p\in\mathbb{R}^{n+1}$ there exists an $H$-bubble passing through $p$. We construct families of nonconstant mappings $H\colon\mathbb{R}^{n+1}\to\mathbb{R}$ of class $C^{1}$ and with further optional properties (like periodicity or prescribed asymptotics at infinity), for which the desired filling property holds. In fact, we use a kind of inverse approach, starting with the construction of suitable bubbles before we exhibit the corresponding functions $H$.

This setup might have some importance in applications to weighted isoperimetric inequalities. For example, looking for area minimizing, disk-type surfaces with prescribed boundary and volume (enclosed within another fixed reference surface with the same boundary), a certain sphere-attaching mechanism is involved. Such construction (due to Wente \cite[p.\,285 ff.]{We1971}, see also \cite{BrCo1984}) rests on the existence of spheres passing through an \emph{a priori} unknown point of the space. The analogous issue in case of prescribed weighted volume would be quite interesting but no result similar to \cite{We1971} seems to be still available. Since prescribing weighted volume is equivalent to prescribing mean curvature, up to a Lagrange multiplier, the presence (or not) of $H$-bubbles filling the space might have some impact with respect to the above described issue.

Let us spend a few words on the main idea developed in this paper. We firstly construct a smooth Jordan curve $\Gamma$ lying in the half-plane of $\mathbb{R}^{n+1}$ defined by $x_{1}=...=x_{n-1}=0$ and $x_{n+1}\ge 0$, passing through the origin, and symmetric about the $x_{n+1}$-axis. Then we introduce a hypersurface of revolution $S$ obtained by revolving the curve $\Gamma$ around the last coordinate axis. Letting $R=\max_{p\in S}|p|$, we define a radial, positive mapping $H=H(r)$ on the ball $B_{R}=\{p\in\mathbb{R}^{n+1}~|~|p|\le R\}$ taking the value of the mean curvature of $S$ at a point $p\in S$ with $|p|=r$. Since $H$ is radially symmetric, every hypersurface obtained by rotation of $S$ about the origin is an $H$-bubble. Hence the $(n+1)$-dimensional ball $B_{R}$ has the filling property with respect to $H$, that is, for every $p\in B_{R}$ there is an $H$-bubble passing through $H$. 

The main difficulty is to construct $\Gamma$ in such a way that $H$ is well defined, of class $C^{1}$, non-constant and with null derivative at $r=0$ and $r=R$. This is achieved by taking $\Gamma$ as a perturbation of a circle. Consequently, the curvature function $H$ turns out to be a perturbation of a positive constant. We point out that the smallness of the perturbation can be suitably controlled (see Example \ref{Example}). 

Since $H$ has null derivative on the boundary of $B_{R}$, we can extend it outside $B_{R}$ in a $C^{1}$ way with the constant value $H(R)$. Then the filling property is automatically satisfied on the whole Euclidean space, because the complement of the ball $B_{R}$ can be filled by round hyperspheres with radius $H(R)^{-1}$. In fact, we can use the function $H$ on $B_{R'}$, for some $R'>R$, as a fundamental block, and arrange infinitely many similar blocks, suitably spaced, to cover the whole space. Also in this case the resulting function $H$ has the property that every point of $\mathbb{R}^{n+1}$ is touched by an $H$-bubble. The arrangement of the blocks or even the form of $H$ at every block can be arbitrarily adjusted in order to fulfill additional requirements on $H$ (like periodicity). Hence, a very huge amount of examples, even with $C^{\infty}$ regularity, can be built. 

The main results are Theorems \ref{T:radial} and \ref{T:main} and are displayed in the last Section, together with an explicit example. The first two sections are devoted to the construction of the profile curve $\Gamma$ and of the reference hypersurface $S$. 

Lastly, we notice that whereas, on one hand, the problem of existence of $H$-bubbles under  global or perturbative conditions on the prescribed mean curvature function requires rather sophisticated tools and arguments (see, e.g., \cite{Tr1985}, \cite{CaMu2011} and the references therein), on the other hand, our inverse approach can be carried out by means of quite elementary methods and one needs just a basic knowledge on surface theory. 

\section{The profile curve}\label{S:profile}
By Jordan curve we mean a closed, simple curve. Here we aim to construct smooth (i.e. $C^{\infty}$) Jordan curves $\Gamma$ in the plane, with the following properties:
\begin{subequations}
\label{Gamma-properties}
\begin{align}
\label{Gamma-symmetry}
&\text{$\Gamma$ is symmetric with respect to the vertical axis}\\
\label{Gamma+}
&0\in\Gamma\quad\text{and}\quad\Gamma\subset\{(x,y)\in\mathbb{R}^{2}~|~y\ge 0\}\\
\label{k-positive}
&\text{$\Gamma$ has positive curvature at every point}\\
\label{|Gamma|}
&\text{for every } r\in[0,R]\text{ there exists a unique }z=(x,y)\in\Gamma\text{ with $x\ge 0$ and $|z|=r$}
\end{align}
\end{subequations}
where $R=\max_{z\in\Gamma}|z|$. More precisely, we prove the following result.
\begin{Lemma}
\label{L:profile}
For every even and $2\pi$-periodic mapping $h\in C^{\infty}(\mathbb{R})$ there exists $\varepsilon_{0}>0$ such that for every $\varepsilon\in(-\varepsilon_{0},\varepsilon_{0})$ the parametric curve
\begin{equation}
\label{Gamma-definition}
\Gamma_{\varepsilon}=\{(g_{\varepsilon}(t)\sin t,g_{\varepsilon}(\pi)+g_{\varepsilon}(t)\cos t)~|~t\in[-\pi,\pi]\,\}\quad\text{with}\quad g_{\varepsilon}:=1+\varepsilon h
\end{equation}
satisfies \eqref{Gamma-properties}. In addition, setting 
\begin{equation}
\label{gamma}
\gamma_{\varepsilon}(t):=(g_{\varepsilon}(t)\sin t,g_{\varepsilon}(\pi)+g_{\varepsilon}(t)\cos t)\,,
\end{equation}
the mapping $|\gamma_{\varepsilon}|$ is strictly decreasing and of class $C^{1}$ in $[0,\pi]$.
\end{Lemma}

\Proof
Fix $h\in C^{\infty}(\mathbb{R})$ even and $2\pi$-periodic and consider the function $g_{\varepsilon}$ and the parametric curve $\Gamma_{\varepsilon}$ as in \eqref{Gamma-definition}. For every $\varepsilon\in\mathbb{R}$ with $|\varepsilon|<\|h\|_{\infty}^{-1}$ we have that $g_{\varepsilon}:=1+\varepsilon h>0$ everywhere and the parametric curve $\Gamma_{\varepsilon}$ defined as in \eqref{Gamma-definition} is a smooth Jordan curve satisfying \eqref{Gamma-symmetry}, with $0\in\Gamma_{\varepsilon}$. Let us show that $\Gamma_{\varepsilon}\subset\{(x,y)\in\mathbb{R}^{2}~|~y\ge 0\}$ if $|\varepsilon|$ is sufficiently small. It is enough to prove that
\begin{equation}
\label{g+}
\tilde{g}_{\varepsilon}(t):=g_{\varepsilon}(\pi)+g_{\varepsilon}(t)\cos t\ge 0\quad\forall t\in[0,\pi]\,.
\end{equation}
One has that
\[
\tilde{g}''_{\varepsilon}(t)=-\cos t+\varepsilon [h(t)\cos t]''\,.
\]
In particular there exists $\varepsilon_{0}>0$ such that $\tilde{g}_{\varepsilon}$ is strictly convex in $[\frac{2\pi}3,\pi]$. Since $g_{\varepsilon}$ is even and $2\pi$-periodic, $g'_{\varepsilon}(\pi)=0$ and then also $\tilde{g}'_{\varepsilon}(\pi)=0$. Therefore $\tilde{g}_{\varepsilon}(t)>\tilde{g}_{\varepsilon}(\pi)=0$ for $t\in[\frac{2\pi}3,\pi]$. In addition
\[
\tilde{g}_{\varepsilon}(t)=-\int_{t}^{\pi}\tilde{g}''_{\varepsilon}(s)(t-s)\,ds=1+\cos t-\varepsilon\int_{t}^{\pi}[h(s)\cos s]''(t-s)\,ds\ge \frac12-\varepsilon C\quad\forall t\in\left[0,\frac{2\pi}3\right]
\]
for some constant $C$ independent of $\varepsilon$ and $t$. Hence, taking $|\varepsilon|$ small, we obtain $\tilde{g}_{\varepsilon}(t)>0$ also for $t\in\left[0,\frac{2\pi}3\right]$. Hence \eqref{g+} holds. 

Now let us consider the parameterization $\gamma_{\varepsilon}$ defined in \eqref{gamma} and let us study the mapping 
\begin{equation}
\label{|gamma|}
|\gamma_{\varepsilon}|(t)=\sqrt{g_{\varepsilon}(t)^{2}+g_{\varepsilon}(\pi)^{2}+2g_{\varepsilon}(\pi)g_{\varepsilon}(t)\cos t}\,.
\end{equation}
We have
\begin{equation}
\label{|gamma|'}
|\gamma_{\varepsilon}|'(t)=\frac{g_{\varepsilon}(t)g'_{\varepsilon}(t)+g_{\varepsilon}(\pi)g'_{\varepsilon}(t)\cos t
-g_{\varepsilon}(\pi)g_{\varepsilon}(t)\sin t
}{\sqrt{g_{\varepsilon}(t)^{2}+g_{\varepsilon}(\pi)^{2}+2g_{\varepsilon}(\pi)g_{\varepsilon}(t)\cos t}}\quad\forall t\in[0,\pi)\,.
\end{equation}
At $t=\pi$ one has to be careful since $|\gamma_{\varepsilon}|(\pi)=0$. Writing the Taylor expansion about $\pi$ yields
\[
|\gamma_{\varepsilon}|'(t)=\frac{g_{\varepsilon}(\pi)^{2}(t-\pi)+o(t-\pi)}{g_{\varepsilon}(\pi)|t-\pi|+o(t-\pi)}\quad\text{as $t\to\pi$.}
\]
Hence $|\gamma_{\varepsilon}|$ is of class $C^{1}$ in $[0,\pi]$. In particular the left derivative at $\pi$ is 
\begin{equation}
\label{|gamma|'pi}
|\gamma_{\varepsilon}|'_{-}(\pi)=-g_{\varepsilon}(\pi)<0\,.
\end{equation}
We also need to control the behaviour of the derivative of $|\gamma_{\varepsilon}|$ at $0$. Writing the Taylor expansion about $0$ yields
\begin{equation}
\label{|gamma|'0}
|\gamma_{\varepsilon}|'(t)=\frac{g_{\varepsilon}(0)g''_{\varepsilon}(0)+g_{\varepsilon}(\pi)g''_{\varepsilon}(0)-g_{\varepsilon}(0)g_{\varepsilon}(\pi)}{g_{\varepsilon}(0)+g_{\varepsilon}(\pi)}\,t+o(t)\quad\text{as~}t\to 0\,.
\end{equation}
In order that the mapping $|\gamma_{\varepsilon}|$ is strictly decreasing we need $|\gamma_{\varepsilon}|_{+}'(0)=0_{-}$. This occurs if 
\begin{equation}
\label{|gamma|'0-}
g_{\varepsilon}(0)g''_{\varepsilon}(0)+g_{\varepsilon}(\pi)g''_{\varepsilon}(0)-g_{\varepsilon}(0)g_{\varepsilon}(\pi)<0\,.
\end{equation} 
Since $g_{\varepsilon}=1+\varepsilon h$, the left-hand side is $-1+O(\varepsilon)$, and then \eqref{|gamma|'0-} holds true taking $\varepsilon>0$ small enough. Hence $|\gamma_{\varepsilon}|'<0$ in a right neighborhood of $0$ and taking a smaller $\varepsilon$ if necessary, in view of \eqref{|gamma|'}--\eqref{|gamma|'pi}, we obtain  the desired monotonicity property for the mapping $|\gamma_{\varepsilon}|$. This immediately implies \eqref{|Gamma|}. 

Finally, let us evaluate the curvature of $\Gamma_{\varepsilon}$. This can be computed by
\begin{equation}
\label{k}
K(\gamma_{\varepsilon}(t))=\frac{i\gamma_{\varepsilon}'(t)\cdot \gamma_{\varepsilon}''(t)}{|\gamma_{\varepsilon}'(t)|^{3}}=\frac{2g_{\varepsilon}'(t)^{2}-g_{\varepsilon}(t)g_{\varepsilon}''(t)+g_{\varepsilon}(t)^{2}}{[(g_{\varepsilon}'(t))^{2}+g_{\varepsilon}(t)^{2}]^{3/2}}\,,
\end{equation}
where, in general, for $z\in\mathbb{R}^{2}$, $iz$ denotes the anticlockwise rotation of $z$ through an angle $\frac\pi2$. Hence, for $g_{\varepsilon}=1+\varepsilon h$, \eqref{k-positive} is fulfilled by taking $|\varepsilon|$ small, since the leading term is $g_{\varepsilon}^{2}$.
\qed  

For future convenience, let us study the regularity property of the curvature of $\Gamma_{\varepsilon}$ as a function of the distance. More precisely, setting 
\begin{equation}
\label{R-epsilon}
R_{\varepsilon}:=\max_{z\in\Gamma_{\varepsilon}}|z|=g_{\varepsilon}(0)+g_{\varepsilon}(\pi)\,,
\end{equation}
let us introduce the mapping ${k}_{\varepsilon}\colon[0,R_{\varepsilon}]\to\mathbb{R}$ defined as 
\begin{equation}
\label{k-epsilon}
{k}_{\varepsilon}:=K\circ\gamma_{\varepsilon}\circ|\gamma_{\varepsilon}|^{-1}
\end{equation}
where $|\gamma_{\varepsilon}|^{-1}$ is the inverse of $|\gamma_{\varepsilon}|\colon[0,\pi]\to[0,R_{\varepsilon}]$. Notice that $|\gamma_{\varepsilon}|^{-1}$ is well defined thanks to Lemma \ref{L:profile}.

\begin{Lemma}
\label{L:curvature}
If $h\in C^{\infty}(\mathbb{R})$ is even, $2\pi$-periodic and satisfies
\begin{equation}
\label{h-derivatives}
h''(0)=h''''(0)=0\,,
\end{equation}
then 
the mapping ${k}_{\varepsilon}$ defined in \eqref{k-epsilon} is of class $C^{1}$ in $[0,R_{\varepsilon}]$ and ${k}_{\varepsilon}'(0)={k}_{\varepsilon}'(R_{\varepsilon})=0$. Moreover, if $h$ vanishes in neighborhoods of $0$ and $\pi$ then $R_{\varepsilon}=2$, ${k}_{\varepsilon}\in C^{\infty}([0,2])$ and $k_{\varepsilon}=1$ in neighborhoods of $0$ and  $2$. 
\end{Lemma}

\Proof
Since $g_{\varepsilon}\in C^{\infty}(\mathbb{R})$, \eqref{k} implies that also $K\circ\gamma_{\varepsilon}\in C^{\infty}(\mathbb{R})$. Moreover $|\gamma_{\varepsilon}|^{-1}$ is of class $C^{1}$ in $[0,R_{\varepsilon})$. Hence also $k_{\varepsilon}$ is so. Moreover, considering that
\begin{equation}
\label{derivative}
\frac{d}{dt}[K\circ\gamma_{\varepsilon}]=
\frac{-3(g_{\varepsilon}')^{3}g_{\varepsilon}''-g_{\varepsilon}(g_{\varepsilon}')^{2}g_{\varepsilon}'''-4g_{\varepsilon}(g_{\varepsilon}')^{3}+3g_{\varepsilon}^{2}g_{\varepsilon}'g_{\varepsilon}''-g_{\varepsilon}^{3}g_{\varepsilon}'''-g_{\varepsilon}^{3}g_{\varepsilon}'+3g_{\varepsilon}g_{\varepsilon}'(g_{\varepsilon}'')^{2}}{[(g_{\varepsilon}')^{2}+g_{\varepsilon}^{2}]^{5/2}}\,,
\end{equation}
and since $|\gamma_{\varepsilon}|^{-1}(0)=\pi$ and $g_{\varepsilon}'(\pi)=g_{\varepsilon}'''(\pi)=0$ (because $g_{\varepsilon}$ is even and $2\pi$-periodic), one obtains $\frac{dk_{\varepsilon}}{dr}(0)=0$. 
Some care is needed at $R_{\varepsilon}$ because $|\gamma_{\varepsilon}|^{-1}(R_{\varepsilon})=0$, $|\gamma_{\varepsilon}|'(0)=0$, whence $\frac{d}{dr}|\gamma_{\varepsilon}|^{-1}(R_{\varepsilon})=\infty$. One has
\[
\frac{dk_{\varepsilon}}{dr}(r)=
\frac{\frac{d(K\circ\gamma_{\varepsilon})}{dt}(|\gamma_{\varepsilon}|^{-1}(r))}{\frac{d|\gamma_{\varepsilon}|}{dt}(|\gamma_{\varepsilon}|^{-1}(r))}\,.
\]
Taking \eqref{|gamma|'0} into account, we need that
\begin{equation}
\label{deriv}
\frac{d(K\circ\gamma_{\varepsilon})}{dt}(t)=o(t)\quad\text{as $t\to 0$}.
\end{equation}
Since $g'_{\varepsilon}(t)=O(t)$, by \eqref{derivative}, equation \eqref{deriv} holds true if
\begin{equation}
\label{4-derivative}
3g_{\varepsilon}(t)^{2}g_{\varepsilon}'(t)g_{\varepsilon}''(t)-g_{\varepsilon}(t)^{3}g_{\varepsilon}(t)'''-g_{\varepsilon}(t)^{3}g_{\varepsilon}'(t)+3g_{\varepsilon}(t)g_{\varepsilon}'(t)g_{\varepsilon}''(t)^{2}=o(t)\quad\text{as $t\to 0$}
\end{equation}
For $g_{\varepsilon}'(0)=g_{\varepsilon}'''(0)=0$, \eqref{4-derivative} is fulfilled when $g_{\varepsilon}''(0)=g_{\varepsilon}''''(0)=0$, too, which follows from \eqref{h-derivatives}. Concerning the last statement, if $h$ vanishes in neighborhoods of $0$ and $\pi$ then $g_{\varepsilon}$ takes the constant value 1 in the same sets. The same holds for $K\circ\gamma_{\varepsilon}$, by \eqref{k}. Hence $\frac{dk_{\varepsilon}}{dr}\in C^{\infty}_{c}((0,R_{\varepsilon}))$, $R_{\varepsilon}=2$ and the assertion follows.
\qed

\section{The reference hypersurface}
\label{S:surface}

In this section we introduce a family of hypersurfaces of revolution $S_{\varepsilon}$ in $\mathbb{R}^{n+1}$ obtained by revolving around the last coordinate axis the curves $\Gamma_{\varepsilon}$ built in Section \ref{S:profile} and placed on the plane $x_{1}=...=x_{n-1}=0$. The symmetry of $\Gamma_{\varepsilon}$ with respect to the vertical axis guarantees that $S_{\varepsilon}$ is well-defined. 

Then we evaluate the principal curvatures $K_{i}$ of $S_{\varepsilon}$ ($i=1,...,n$) and we impose conditions on the parameterization of $\Gamma_{\varepsilon}$ which guarantee that the functions $K_{i}$ depend in a $C^{1}$ way on the distance from the origin, with null derivative at the origin and at the maximum distance. 

To carry out this plan, it is convenient to parametrize the $n$-dimensional unit sphere $\mathbb{S}^{n}$ in the Euclidean $(n+1)$-space in hyperspherical coordinates, as follows:
\[
\sigma(\theta_{1},...,\theta_{n})=\left[\begin{array}{c}
\sin\theta_{1}\,\sin\theta_{2}\cdots\,\sin\theta_{n}\\
\cos\theta_{1}\,\sin\theta_{2}\cdots\,\sin\theta_{n}\\
\cos\theta_{2}\,\sin\theta_{3}\cdots\,\sin\theta_{n}\\
\vdots\\
\cos\theta_{n-1}\,\sin\theta_{n}\,\\
\cos\theta_{n}\end{array}\right]
\]
with $\theta_{1}\in[0,2\pi]$, $\theta_{i}\in[0,\pi]$ for $i=2,...,n$. 

Fixing a mapping $g\colon\mathbb{R}\to(0,\infty)$ of class $C^{\infty}$, even and $2\pi$-periodic, let us introduce the hypersurface $S_{g}$ parameterized by
\[
x(\theta_{1},...,\theta_{n})=g(\theta_{n})\sigma(\theta_{1},...,\theta_{n})\,,\quad(\theta_{1},...,\theta_{n})\in[0,2\pi]\times[0,\pi]^{n-1}\,.
\]

\begin{Lemma}
\label{L:Ki}
For every $p=x(\theta_{1},...,\theta_{n})\in S_{g}$ the principal curvatures of $S_{g}$ at $p$ are given by
\begin{equation}
\label{principal-curvatures}
K_{1}=...=K_{n-1}=\frac{g-\frac{\cos\theta_{n}}{\sin\theta_{n}}\,g'}{g\sqrt{g^{2}+(g')^{2}}}\quad\text{and}\quad K_{n}=\frac{g^{2}+2(g')^{2}-gg''}{[g^{2}+(g')^{2}]^{3/2}}
\end{equation}
where the primes denote derivatives with respect to $\theta_{n}$.
\end{Lemma}

\Proof
Let us fix a point $p=x(\theta_{1},...,\theta_{n})\in S_{g}$. For the sake of brevity, we suppress the variables $\theta_1,...,\theta_n$ in our notation. Let us introduce the vectors
\[
\tau_{i}=\frac{\partial\sigma}{\partial\theta_{i}}\quad(i=1,...,n).
\]
One can check that 
\begin{equation}
\label{orthogonal}
\sigma\cdot\tau_{i}=\tau_{i}\cdot\tau_{j}=0\quad\forall i,j=1,...,n,~i\ne j\,.
\end{equation}
Thanks to \eqref{orthogonal}, the outward-pointing normal versor $N$ at $p$ can be expressed in the form
\[
N=\alpha_{0}\sigma+\sum_{i=1}^{n}\alpha_{i}\tau_{i}
\]
where $\alpha_{0},...,\alpha_{n}\in\mathbb{R}$ satisfy
\begin{equation}
\label{Ort}
N\cdot\frac{\partial x}{\partial \theta_{i}}=0\quad\forall i=1,...,n\,.
\end{equation}
We have
\[
\begin{split}
&N\cdot\frac{\partial x}{\partial \theta_{i}}=\left(\alpha_{0}\sigma+\sum_{j=1}^{n}\alpha_{j}\tau_{j}\right)\cdot g\tau_{i}=\alpha_{i}g|\tau_{i}|^{2}\qquad(i=1,...,n-1)
\\
&N\cdot\frac{\partial x}{\partial \theta_{n}}=\left(\alpha_{0}\sigma+\sum_{j=1}^{n}\alpha_{j}\tau_{j}\right)\cdot (g'\sigma+g\tau_{n})=\alpha_{0}g'+\alpha_{n}g
\end{split}
\]
since $|\sigma|=|\tau_{n}|=1$. Hence the equations \eqref{Ort} imply
\[
\alpha_{1}=\cdots=\alpha_{n-1}=0\,,\quad
\alpha_{n}=-\frac{\alpha_{0}g'}{g}
\]
and then we can write the outward-pointing normal versor in the form
\[
N=\frac{1}{\sqrt{g^{2}+(g')^{2}}}\left(g\sigma-g'\tau_{n}\right)\,.
\]
We compute
\[
\begin{split}
&\frac{\partial(N\circ x)}{\partial\theta_{i}}=\frac{g}{\sqrt{g^{2}+(g')^{2}}}\tau_{i}-\frac{g'}{\sqrt{g^{2}+(g')^{2}}}\frac{\partial\tau_{n}}{\partial\theta_{i}}\qquad(i=1,...,n-1)
\\
&\frac{\partial(N\circ x)}{\partial\theta_{n}}=\frac{g}{\sqrt{g^{2}+(g')^{2}}}\tau_{n}-\frac{g'}{\sqrt{g^{2}+(g')^{2}}}\frac{\partial\tau_{n}}{\partial\theta_{n}}+\left[\frac{g}{\sqrt{g^{2}+(g')^{2}}}\right]'\sigma-\left[\frac{g'}{\sqrt{g^{2}+(g')^{2}}}\right]'\tau_{n}\,.
\end{split}
\]
Since 
\[
\frac{\partial\tau_{n}}{\partial\theta_{i}}=\frac{\cos\theta_{n}}{\sin\theta_{n}}\,\tau_{i}\quad(i=1,...,n-1)\quad\text{and}\quad\frac{\partial\tau_{n}}{\partial\theta_{n}}=-\sigma\,,
\]
we obtain
\begin{equation}
\label{ki}
\begin{split}
&\frac{\partial(N\circ x)}{\partial\theta_{i}}=\frac{g-\frac{\cos\theta_{n}}{\sin\theta_{n}}\,g'}{g\sqrt{g^{2}+(g')^{2}}}\,\frac{\partial x}{\partial\theta_{i}}\quad(i=1,...,n-1)
\\
&\frac{\partial(N\circ x)}{\partial\theta_{n}}=\frac{g^{2}+2(g')^{2}-gg''}{[g^{2}+(g')^{2}]^{3/2}}\,\frac{\partial x}{\partial\theta_{n}}\,.
\end{split}
\end{equation}
By definition, the principal curvatures at $p\in S_{g}$ are the eigenvalues of the shape operator $L_{p}\colon T_{p}S_{g}\to T_{p}S_{g}$ given by
\[
v\mapsto L_{p}v=\frac{\partial N}{\partial v}(p)\quad(v\in T_{p}S_{g})\,.
\]
Since, for $p=x(\theta_{1},...,\theta_{n})$ and for every $i=1,...,n$ one has 
\[
L_{p}\frac{\partial x}{\partial\theta_{i}}=\frac{\partial(N\circ x)}{\partial\theta_{i}}\,,
\]
by \eqref{ki} the principal curvatures of $S$ are given by \eqref{principal-curvatures}.
\qed

Fixing $h\in C^{\infty}(\mathbb{R})$ even and $2\pi$-periodic, for $\varepsilon\in\mathbb{R}$ let $g_{\varepsilon}$ be as in \eqref{Gamma-definition} and let $S_{\varepsilon}$ be the $n$-dimensional hypersurface in $\mathbb{R}^{n+1}$ parametrized by
\[
x_{\varepsilon}(\theta_{1},...,\theta_{n})=g_{\varepsilon}(\theta_{n})\sigma(\theta_{1},...,\theta_{n})+g_{\varepsilon}(\pi)\mathbf{e}_{n+1}\quad(\theta_{1},...,\theta_{n})\in[0,2\pi]\times[0,\pi]^{n-1}
\]
where $\mathbf{e}_{n+1}=\sigma(0,...,0)$. 
Observe that $S_{\varepsilon}$ is obtained from $S_{g_\epsilon}$ via translation by an amount of $g_\epsilon(\pi)$ in the direction of the last coordinate axis and satisfies the following properties.

\begin{Lemma}
\label{L:surface}
There exists $\varepsilon_{0}>0$ such that for every $\varepsilon\in(-\varepsilon_{0},\varepsilon_{0})$:
\begin{itemize}
\item[(i)]
$S_{\varepsilon}$ is a compact, embedded $n$-dimensional hypersurface of revolution around the last coordinate axis of $\mathbb{R}^{n+1}$, diffeomorphic to $\mathbb{S}^{n}$;
\item[(ii)]
$0\in S_{\varepsilon}$ and $S_{\varepsilon}$ is contained in the half-space $\{(x_{1},...,x_{n+1})\in\mathbb{R}^{n+1}~|~x_{n+1}\ge 0\}$;
\item[(iii)]
$\max_{p\in S_{\varepsilon}}|p|=R_{\varepsilon}$ where $R_{\varepsilon}>0$ is given by \eqref{R-epsilon};
\item[(iv)]
the principal curvatures of $S_{\varepsilon}$ satisfy $K_{i}(p)>0$ for all $p\in S_{\varepsilon}\setminus\{0,R_{\varepsilon}\mathbf{e}_{n+1}\}$ ($i=1,...,n-1$) and $K_{n}(p)>0$ for all $p\in S_{\varepsilon}$.
\end{itemize}
\end{Lemma}

\Proof
Properties \emph{(i)--(iii)} follow from the definition of $S_{\varepsilon}$ and from Lemma \ref{L:profile}. Let us discuss \emph{(iv)}. According to \eqref{k} and \eqref{principal-curvatures}, the $n$-th principal curvature $K_{n}$ of $S_{\varepsilon}$ equals the curvature of $\Gamma_{\varepsilon}$. Therefore $K_{n}>0$ on $S_{\varepsilon}$, again by Lemma \ref{L:profile}. By \eqref{principal-curvatures}, for $i=1,...,n-1$, one has $K_{i}>0$ on $S_{\varepsilon}\setminus\{0,R_{\varepsilon}\mathbf{e}_{n+1}\}$ when 
\begin{equation}
\label{Ki-positive}
g_{\varepsilon}(t)\sin t-g'_{\varepsilon}(t)\cos t>0\quad\forall t\in(0,\pi).
\end{equation} 
One has 
\[
g_{\varepsilon}(t)\sin t-g'_{\varepsilon}(t)\cos t=\left[1+\varepsilon\widetilde{h}(t)\right]\sin t\quad\text{where }\widetilde{h}(t)=h(t)-\frac{h'(t)}{\sin t}\cos t\,.
\]
Since $h'(0)=h'(\pi)=0$,
\[
\lim_{t\to 0}\frac{h'(t)}{\sin t}=h''(0)\quad\text{and}\quad\lim_{t\to\pi}\frac{h'(t)}{\sin t}=h''(\pi)\,,
\]
the mapping $\widetilde{h}$ is continuous in $[0,\pi]$. Therefore, for $|\varepsilon|$ sufficiently small, \eqref{Ki-positive} holds true.
\qed

Let us study the regularity property of the principal curvatures of $S_{\varepsilon}$ as functions of the distance from the origin. To this aim, let us introduce the mapping $y_{\varepsilon}\colon[0,R_{\varepsilon}]\to S_{\varepsilon}$ defined by
\begin{equation}
\label{y-eps}
y_{\varepsilon}(r):=x_{\varepsilon}(0,...,0,|\gamma_{\varepsilon}|^{-1}(r))
=\left[g_{\varepsilon}(t)(\sin t)\mathbf{e}_{n}+(g_{\varepsilon}(t)\cos t+g_{\varepsilon}(\pi))\mathbf{e}_{n+1}\right]_{t=|\gamma_{\varepsilon}|^{-1}(r)}\,,
\end{equation}
where $\gamma_{\varepsilon}$ is defined in \eqref{gamma} and $|\gamma_{\varepsilon}|^{-1}$ is the inverse of the mapping $|\gamma_{\varepsilon}|\colon[0,\pi]\to[0,R_{\varepsilon}]$ (see Lemma \ref{L:profile}). 
\begin{Lemma}
\label{L:Ki-regular}
If \eqref{h-derivatives} and 
\begin{equation}
\label{h-derivative-pi}
h''(\pi)=0
\end{equation}
hold, then 
the mappings $K_{i}\circ y_{\varepsilon}\colon[0,R_{\varepsilon}]\to\mathbb{R}$ ($i=1,...,n$) are of class $C^{1}$ in $[0,R_{\varepsilon}]$ and with null derivatives at end points. Moreover, if $h$ vanishes in neighborhoods of $0$ and $\pi$ then $R_{\varepsilon}=2$, $K_{i}\circ y_{\varepsilon}\in C^{\infty}([0,2])$ and $K_{i}\circ y_{\varepsilon}=1$ in neighborhoods of $0$ and $2$ ($i=1,...,n$). 
\end{Lemma}

\Proof We observe that the mapping $K_{n}\circ y_{\varepsilon}$ is the function $k_{\varepsilon}$ defined in \eqref{k-epsilon} and then the regularity of $K_{n}\circ y_{\varepsilon}$ immediately follows from Lemma \ref{L:curvature}. Let us consider $K_{i}\circ y_{\varepsilon}$ with $i=1,...,n-1$. Setting 
\begin{equation}
\label{bar-k-eps}
\bar{k}(t):=\frac{g_{\varepsilon}(t)\sin t-g'_{\varepsilon}(t)\cos t}{g_{\varepsilon}(t)\sin t\sqrt{g_{\varepsilon}(t)^{2}+g'_{\varepsilon}(t)^{2}}}
\end{equation}
we have that $K_{i}\circ y_{\varepsilon}=\bar{k}\circ|\gamma_{\varepsilon}|^{-1}$ ($i=1,...,n-1$). Moreover $\bar{k}\in C^{1}((0,\pi))$. Therefore $K_{i}\circ y_{\varepsilon}\in C^{1}((0,R_{\varepsilon}))$, with
\[
\frac{d[K_{i}\circ y_{\varepsilon}]}{dr}(r)=\frac{\bar{k}'(|\gamma_{\varepsilon}|^{-1}(r))}{|\gamma_{\varepsilon}|'(|\gamma_{\varepsilon}|^{-1}(r))}\,.
\]
We aim to show that $\frac{d[K_{i}\circ y_{\varepsilon}]}{dr}(r)\to 0$ as $r\to 0$ and as $r\to R_{\varepsilon}$. Let us study the limit as $r\to 0$. We have
\[
\lim_{r\to 0}\frac{d[K_{i}\circ y_{\varepsilon}]}{dr}(r)=\lim_{t\to\pi}\frac{\bar{k}'(t)}{|\gamma_{\varepsilon}|'(t)}\,.
\]
Moreover
\[\begin{split}
\bar{k}'(t)
&=
\left[\frac{g_{\varepsilon}'(t)}{\sin t}\right]^{2}\!\!\!\frac{\sin t\,\cos t}{g_{\varepsilon}(t)^{2}\sqrt{g_{\varepsilon}(t)^{2}+g_{\varepsilon}'(t)^{2}}}+\frac{g_{\varepsilon}'(t)-g_{\varepsilon}''(t)\,\sin t\,\cos t}{(\sin t)^{2}\,g_{\varepsilon}(t)\sqrt{g_{\varepsilon}(t)^{2}+g_{\varepsilon}'(t)^{2}}}\\
&\quad-\frac{g_{\varepsilon}'(t)[g_{\varepsilon}(t)+g_{\varepsilon}''(t)]}{(\sin t)\,g_{\varepsilon}(t)[g_{\varepsilon}(t)^{2}+g_{\varepsilon}'(t)^{2}]^{3/2}}\,.
\end{split}
\]
Since $g_{\varepsilon}'(\pi)=0$, we have that $g_{\varepsilon}'(t)-g_{\varepsilon}''(t)\,\sin t\,\cos t=O((t-\pi)^{3})$ and, by the rule of de l'H\^opital, $\frac{g_{\varepsilon}'(t)}{\sin t}\to -g_{\varepsilon}''(\pi)$ as $t\to\pi$.
In addition $|\gamma_{\varepsilon}|'(t)\to-g_{\varepsilon}(\pi)<0$ as $t\to\pi_{-}$ (see \eqref{|gamma|'pi}). Therefore
\[
\lim_{r\to 0}\frac{d[K_{i}\circ y_{\varepsilon}]}{dr}(r)=-\frac{g_{\varepsilon}''(\pi)[g_{\varepsilon}(\pi)+g_{\varepsilon}''(\pi)]}{g_{\varepsilon}(\pi)^{4}}=0
\]
because $g_{\varepsilon}''(\pi)=\varepsilon h''(\pi)=0$ by \eqref{h-derivative-pi}. 

Now let us study the limit as $r\to R_{\varepsilon}$. We have
\[
\lim_{r\to R_{\varepsilon}}\frac{d[K_{i}\circ y_{\varepsilon}]}{dr}(r)=\lim_{t\to 0}\frac{\tilde{K}_{1}'(t)}{|\gamma_{\varepsilon}|'(t)}\,.
\]
We know that
\[
|\gamma_{\varepsilon}|'(t)=\frac{[g_{\varepsilon}(0)g''_{\varepsilon}(0)+g_{\varepsilon}(\pi)g''_{\varepsilon}(0)-g_{\varepsilon}(0)g_{\varepsilon}(\pi)]}{g_{\varepsilon}(0)+g_{\varepsilon}(\pi)}\,t+o(t)\quad\text{as~}t\to 0
\]
(see \eqref{|gamma|'0}). We also know that $g_{\varepsilon}'(0)=g_{\varepsilon}''(0)=g_{\varepsilon}'''(0)=g_{\varepsilon}''''(0)=0$ beacuse of \eqref{h-derivatives}. Then, by the rule of de l'H\^opital, $\frac{g_{\varepsilon}'(t)}{\sin t}\to -g_{\varepsilon}''(\pi)=0$ as $t\to 0$. Moreover $g_{\varepsilon}'(t)-g_{\varepsilon}''(t)\sin t\,\cos t=o(t^{4})$ as $t\to 0$. Then
\[
\frac{g_{\varepsilon}'(t)-g_{\varepsilon}''(t)\sin t\,\cos t}{(\sin t)^{2}\,g_{\varepsilon}(t)\sqrt{g_{\varepsilon}^{2}(t)+g_{\varepsilon}'(t)^{2}}}=o(t^{2})\quad\text{as }t\to 0\,.
\]
In addition also $g_{\varepsilon}'(t)=o(t^{4})$ as $t\to 0$. Therefore
\[
\frac{g_{\varepsilon}'(t)[g_{\varepsilon}(t)+g_{\varepsilon}''(t)]}{(\sin t)\,g_{\varepsilon}(t)[g_{\varepsilon}(t)^{2}+g_{\varepsilon}'(t)^{2}]^{3/2}}=o(t^{2})\quad\text{as }t\to 0\,.
\]
In conclusion, also 
\[
\lim_{r\to R_{\varepsilon}}\frac{d[K_{i}\circ y_{\varepsilon}]}{dr}(r)=0\,.
\]
For the last part of the Lemma, one argues exactly as in the last part of the proof of Lemma \ref{L:curvature}, taking account of \eqref{principal-curvatures}. 
\qed

\section{Families of curvature functions with the filling property}
In this section we construct families of non-constant mappings $H\colon\mathbb{R}^{n+1}\to\mathbb{R}$ having the following properties:
\begin{itemize}
\item[(H)$_{1}$] $H\in C^{1}(\mathbb{R}^{n+1})$;
\item[(H)$_{2}$] there exist $C_{1},C_{2}>0$ such that $C_{1}\le H(p)\le C_{2}$ for every $p\in\mathbb{R}^{n+1}$;
\item[(H)$_{3}$] for every $p\in\mathbb{R}^{n+1}$ there exists an embedded hypersurface $S$ diffeomorphic to $\mathbb{S}^{n}$ with mean curvature $H$ at every point and with $p\in S$.
\end{itemize}
Other optional properties on $H$ like periodicity or some asymptotic behaviour can be added. 

As a first result, we exhibit a family of radially symmetric curvature functions which satisfy (H)$_{1}$--(H)$_{3}$ and are constant outside a ball. 

\begin{Theorem}
\label{T:radial}
Let $h\colon\mathbb{R}\to\mathbb{R}$ be a $2\pi$-periodic, even function of class $C^{\infty}$ satisfying \eqref{h-derivatives} and \eqref{h-derivative-pi}. 
For $\varepsilon\in\mathbb{R}$ let $g_{\varepsilon}$, $R_{\varepsilon}$, $K_{i}$, and $y_{\varepsilon}$ as in \eqref{Gamma-definition}, \eqref{R-epsilon}, \eqref{principal-curvatures} with $g=g_{\varepsilon}$, and \eqref{y-eps}, respectively. Then there exists $\varepsilon_{0}>0$ such that for every $\varepsilon\in(-\varepsilon_{0},\varepsilon_{0})$ the function $H_{\varepsilon}\colon\mathbb{R}^{n+1}\to\mathbb{R}$ defined by
\begin{equation}
\label{H-eps}
H_{\varepsilon}(p):=
\begin{cases}
\frac{g_{\varepsilon}(0)}{n}\sum_{i=1}^{n}(K_{i}\circ y_{\varepsilon})(|p|)&\text{for $|p|\le R_{\varepsilon}$}\\
1&\text{for $|p|>R_{\varepsilon}$.}
\end{cases}
\end{equation}
satisfies \emph{(H)$_{1}$--(H)$_{3}$} and is radially symmetric. In addition, if $\varepsilon\to 0$ then $R_{\varepsilon}\to 2$ and $H_{\varepsilon}\to 1$ in $C^{1}(\mathbb{R}^{n+1})$. Moreover, if $h$ vanishes in neighborhoods of $0$ and $\pi$ then $H_{\varepsilon}$ is of class $C^{\infty}$.  
\end{Theorem}

\Proof
Fixing $h$ as in the statement of the Theorem, and taking $\varepsilon\in\mathbb{R}$ with $|\varepsilon|$ small enough, according to Lemma \ref{L:surface}, the hypersurfaces ${S}_{\varepsilon}$ built in Section \ref{S:surface} is diffeomorphic to $\mathbb{S}^{n}$ and, for every $(\theta_{1},...,\theta_{n})\in[0,2\pi]\times[0,\pi]^{n-1}$ the mean curvature of $S_{\varepsilon}$ at $p=x_{\varepsilon}(\theta_{1},...,\theta_{n})\in S_{\varepsilon}$ is given by 
\[
M=\frac{1}{n}\sum_{i=1}^{n}K_{i}
\]
with $K_{i}$ as in \eqref{principal-curvatures}. Let
\[
\widetilde{H}_{\varepsilon}(r):=(M\circ y_{\varepsilon})(r)\quad\forall r\in[0,R_{\varepsilon}]
\]
with $y_{\varepsilon}$ as in \eqref{y-eps}. Thus $\widetilde{H}_{\varepsilon}(|p|)$ equals the mean curvature of $S_{\varepsilon}$ at $p$. By radial symmetry, the same holds true for any hypersurface obtained by rotating $S_{\varepsilon}$ about the origin of $\mathbb{R}^{n+1}$. Hence, for every $p\in\mathbb{R}^{n+1}$ with $|p|\le R_{\varepsilon}$ there exists an embedded hypersurface diffeomorphic to $\mathbb{S}^{n}$ with mean curvature $H$ at every point and with $p\in S$. By Lemma \ref{L:Ki-regular}, $\widetilde{H}_{\varepsilon}\in C^{1}([0,R_{\varepsilon}])$ with $\widetilde{H}'_{\varepsilon}(0)=\widetilde{H}'_{\varepsilon}(R_{\varepsilon})=0$ and setting $\widetilde{H}_{\varepsilon}(r):=\widetilde{H}_{\varepsilon}(R_{\varepsilon})$ for $r>R_{\varepsilon}$, we obtain a function of class $C^{1}$ on $[0,\infty)$. In particular $\widetilde{H}_{\varepsilon}(R_{\varepsilon})=\frac{1}{g_{\varepsilon}(0)}$ and, by Lemma \ref{L:surface} (iv), there exist $C_{1},C_{2}>0$ such that $C_{1}\le \widetilde{H}_{\varepsilon}(r)\le C_{2}$ for every $r>0$. Moreover for every point $p\in\mathbb{R}^{n+1}$ with $|p|\le R_{\varepsilon}$ one can take a round hypersphere of radius $g_{\varepsilon}(0)$ whose mean curvature equals $\widetilde{H}_{\varepsilon}$. Finally, because of the above discussion, the function $H_{\varepsilon}(p)=g_{\varepsilon}(0)\widetilde{H}_{\varepsilon}(|p|)$, as defined in \eqref{H-eps}, satisfies {(H)$_{1}$--(H)$_{3}$} and is radially symmetric. The last properties plainly follow, because of the definition of $g_{\varepsilon}$ and by arguing as in the last part of the proof of Lemma \ref{L:curvature}. 
\qed

\begin{Theorem}
\label{T:main}
Let $h\colon\mathbb{R}\to\mathbb{R}$ be a $2\pi$-periodic, even function of class $C^{\infty}$ satisfying \eqref{h-derivatives} and \eqref{h-derivative-pi}. Then there exists $\varepsilon_{0}>0$ such that for every set of numbers $\{\varepsilon_{j}\}_{j\in\mathbb{N}}\subset\mathbb{R}$ with 
\[
|\varepsilon_{j}|<\varepsilon_{0}\quad\forall j\in\mathbb{N}
\] 
and for every set of points $\{p_{j}\}_{j\in\mathbb{N}}\subset\mathbb{R}^{n+1}$ such that 
\begin{equation}
\label{pi-pj}
|p_{i}-p_{j}|\ge \max\{R_{\varepsilon_{i}},R_{\varepsilon_{j}}\}+2
\end{equation}
the function $H\colon\mathbb{R}^{n+1}\to\mathbb{R}$ defined by
\begin{equation}
\label{H-global}
H(p)=\begin{cases}H_{\varepsilon_{j}}(p-p_{j})&\text{ if $|p-p_{j}|\le R_{\varepsilon_{j}}$ for some $j\in\mathbb{N}$}\\
1&\text{otherwise}
\end{cases}
\end{equation}
with $H_{\varepsilon}$ as in \eqref{H-eps} satisfies \emph{(H)$_{1}$--(H)$_{3}$}. In particular, if $\varepsilon_{j}=\varepsilon$ for all $j$ and the set $\{p_{j}\}$ is periodic (in all directions), then the mapping $H$ is periodic. If $\varepsilon_{j}\to 0$ as $j\to\infty$ then $H(p)\to 1$ as $|p|\to\infty$. Moreover if $h$ vanishes in neighborhoods of $0$ and $\pi$ then $H$ is of class $C^{\infty}$. 
\end{Theorem}

\Proof The result plainly follows from Theorem \ref{T:radial} and from the fact that the gluing of the blocks defined by $H(x)=H_{\varepsilon_{j}}(x -p_{j})$ on $B_{R_{\varepsilon_{j}}}(p_{j})$ outside the region $\bigcup_{i}B_{R_{\varepsilon_{i}}}(p_{i})$ is nice because this function takes the common value 1 outside $B_{R_{\varepsilon_{j}}}(p_{j})$. The filling property is satisfied because it holds in each ball $B_{R_{\varepsilon_{j}}+2}(p_{j})$ and these balls are pairwise disjoint by \eqref{pi-pj}. Hence the region $\mathbb{R}^{n+1}\setminus \bigcup_{i}B_{R_{\varepsilon_{i}}}(p_{i})$ can be filled by round unit spheres. 
\qed

\begin{Remark}
Even more complicated mappings satisfying {(H)$_{1}$--(H)$_{3}$} can be constructed by taking a sequence of $2\pi$-periodic, even functions $h_{j}\in C^{\infty}(\mathbb{R}\to\mathbb{R})$ ($j\in\mathbb{N}$), considering a corresponding sequence of maps $H_{\varepsilon_{j},h_{j}}$ defined as in \eqref{H-eps} and then gluing them according to \eqref{H-global}. 
\end{Remark}

\begin{Example}
\label{Example}
As a function $h$ which satisfies the assumptions of Theorems \ref{T:radial} and \ref{T:main} one can take $h(t)=\sin^{6}t$. One can also give an estimate on the interval of admissible values for the smallness parameter $\varepsilon$. In fact, considering the proofs of Lemmata of Sections \ref{S:profile} and \ref{S:surface}, one needs $|\varepsilon|<\|h\|_{\infty}^{-1}$ such that:
\begin{itemize} 
\item
$g_{\varepsilon}(\pi)+g_{\varepsilon}(t)\cos t>0$ for $t\in[0,\pi)$ (see \eqref{g+});
\item
the mapping $|\gamma_{\varepsilon}|$ defined in \eqref{|gamma|} is strictly decreasing in $[0,\pi]$;
\item
the curvatures written in \eqref{k} and \eqref{bar-k-eps} are positive, respectively, in $[0,\pi]$ and in $(0,\pi)$.
\end{itemize}
Taking $h(t)=(\sin t)^{6}$, with elementary computations one can check that the previous conditions are fulfilled taking $g_{\varepsilon}=1+\varepsilon h$ with $\varepsilon \in (-\frac17,\frac25)$. 
\end{Example}

\noindent

\end{document}